\documentclass[10pt]{article}

\usepackage[margin=1.3in]{geometry}
\usepackage{amsmath,amssymb}
\usepackage{changepage}
\usepackage[utf8x]{inputenc}
\usepackage{textcomp,marvosym}
\usepackage{cite}
\usepackage{nameref,hyperref}
\usepackage[right]{lineno}
\usepackage{microtype}
\DisableLigatures[f]{encoding = *, family = * }
\usepackage[table]{xcolor}
\usepackage{array}
\newcolumntype{+}{!{\vrule width 2pt}}

\newlength\savedwidth

\usepackage[aboveskip=1pt,labelfont=bf,labelsep=period,justification=raggedright,singlelinecheck=off]{caption}

\makeatletter
\renewcommand{\@biblabel}[1]{\quad#1.}
\makeatother

\usepackage{lastpage,fancyhdr,graphicx}
\usepackage{epstopdf}

\usepackage{graphicx}
\usepackage{color}
\usepackage{subcaption}
\usepackage{hyperref}
\usepackage{paralist}
\usepackage{amsmath}
\usepackage{amsthm}
\usepackage{amsfonts}
\usepackage{amssymb}
\usepackage{algorithm}
\usepackage{algpseudocode}

\theoremstyle{plain}

\theoremstyle{remark}

\theoremstyle{definition}

\usepackage{soul,xcolor}
\setstcolor{red}

\begin{document}
\vspace*{0.2in}

\begin{flushleft}
{\Large
\textbf\newline{Exponentially decaying modes and long-term prediction of sea ice concentration using Koopman Mode Decomposition}
}
\newline
\\
James Hogg\textsuperscript{1},
Maria Fonoberova\textsuperscript{1*},
Igor Mezi\'c\textsuperscript{1,2}
\\
\bigskip
\textbf{1}  Aimdyn, Inc., Santa Barbara, CA 93101, USA
\\
\textbf{2} University of California - Santa Barbara
\\
\bigskip
* mfonoberova@aimdyn.com
\end{flushleft}

\section*{Abstract}

Sea ice cover in the Arctic and Antarctic is an important indicator of changes in the climate, with important environmental, economic and security consequences. The complexity of the spatio-temporal dynamics of sea ice makes it difficult to assess the temporal nature of the changes - e.g. linear or exponential - and their precise geographical loci. In this study, Koopman Mode Decomposition (KMD) was applied to satellite data of sea ice concentration for the northern and southern hemispheres to gain insight into the temporal and spatial dynamics of the sea ice behavior and to predict future sea ice behavior. We discover exponentially decaying spatial modes in both hemispheres and discuss their precise spatial extent, and also perform precise geographic predictions of sea ice concentration up to four years in the future. This data-driven decomposition technique gives insight in spatial and temporal dynamics not apparent in traditional linear approaches.

\section*{Introduction}

Sea ice is floating ice that forms when ocean water freezes. The formation and distribution of sea ice plays an important role in the planet's climate and thus large amounts of data related to quantitative measures of sea ice have been collected, including continuous satellite remote sensing measurements since 1978. The decreasing extent of Arctic sea ice extent over the last several decades has had negative effects on Arctic wildlife and local communities, while also potentially opening new regions to maritime commerce and natural resources exploration. The future of sea ice behavior is thus of great significance for environmental, economic, and national security reasons. There are several studies that suggest a nonlinear trend in the decline of the sea ice cover \cite{comisoetal:2008,stroeveetal:2012}. 

A variety of approaches have been applied to predict future sea ice behavior over short time scales (1-3 months in the future), including both dynamical (model-based, either coupled ice-ocean-atmosphere or ice-ocean models) and statistical (data-based) (see the reports of the Sea Ice Prediction Network (SIPN) \cite{MeierEtAl:2018} and, e.g. kernel analog forecasting \cite{ComeauEtAl:2018}). Statistical based approaches were reported by the SIPN to be similarly or more successful compared with dynamical approaches for prediction of sea ice distributions 1-3 months in the future \cite{MeierEtAl:2018}. Of statistical based approaches, spatial mode based approaches \cite{KondrashovEtAl:2018} have particular advantages. Compared to regression or trend based approaches, spatial mode based approaches are powerful tools for studying and predicting the geographic and temporal behavior of sea ice because they decompose the time dependent sea ice data into time varying spatial structures of physical significance. Here we apply Koopman Mode Decomposition (KMD) \cite{Mezic:2005} to sea ice concentration dynamics and prediction. 

KMD is a mathematical tool well suited to analyzing sea ice dynamical behavior because it identifies important spatial structures and their complex time dependence from large data sets such as those available for sea ice. The Koopman operator theory \cite{MezicandBanaszuk:2004,Mezic:2005,Rowleyetal:2009,Budisicetal:2012,Williamsetal:2015,Bruntonetal:2016,Giannakisetal:2015} is already widely used to analyze data and provide models for complex dynamic processes. Mathematically, the Koopman operator \cite{Nageletal:2014} is a linear representation on an observables space of an appropriate group action on state space. As it is a linear representation, it is natural that key objects of analysis will be its eigenvalues and eigenfunctions. It turns out that in distributed systems (for example, those with a spatial component, e.g. a fluid velocity field, or a dynamical system on a graph), that, in addition to the eigenvalues and eigenfunctions, there is a third class of objects of importance: the Koopman modes \cite{Rowleyetal:2009}. The eigenvalues of the Koopman operator provide the time scales on which - a potentially exponential - change in sea ice cover is happening, and the modes indicate the spatial extent of the changes. Crucially, the Koopman operator methods do not require a model – observables like sea ice thickness and concentration are sufficient to compute the eigenvalues and the associated modes.

The most popular computational method for Koopman eigenvalues and modes is the Dynamic Mode Decomposition (DMD), that has become a major tool in the data-driven analysis of complex dynamical systems. DMD was first introduced in 2008 by P. Schmid \cite{Schmid:2008wv} for the study of fluid flows where it was conceptualized as an algorithm to decompose the flow field into component fluid structures, called “dynamic modes” or “DMD modes”, that described the evolution of the flow. The DMD modes and their temporal behavior are given by the spectral analysis of the linear operator, which is constructed from data since it is assumed that direct access to it is not available. The book \cite{DMDbook} provides references and introduction to a variety of DMD-related algorithms. Rowley et al. \cite{Rowleyetal:2009} gave the method theoretical underpinnings by connecting it to the spectral analysis of the Koopman operator. The paper \cite{Drmacetal:2018} provides both enhancements and analysis to the DMD method as well as additional theoretical underpinning for its relationship to the Koopman operator, and contrast to the Galerkin projection (or finite section) methods such as EDMD \cite{Williamsetal:2015}.

\begin{figure}[h!]
\centering
\begin{subfigure}[b]{0.33\textwidth}
\centering
\fbox{\includegraphics[width=0.95\textwidth]{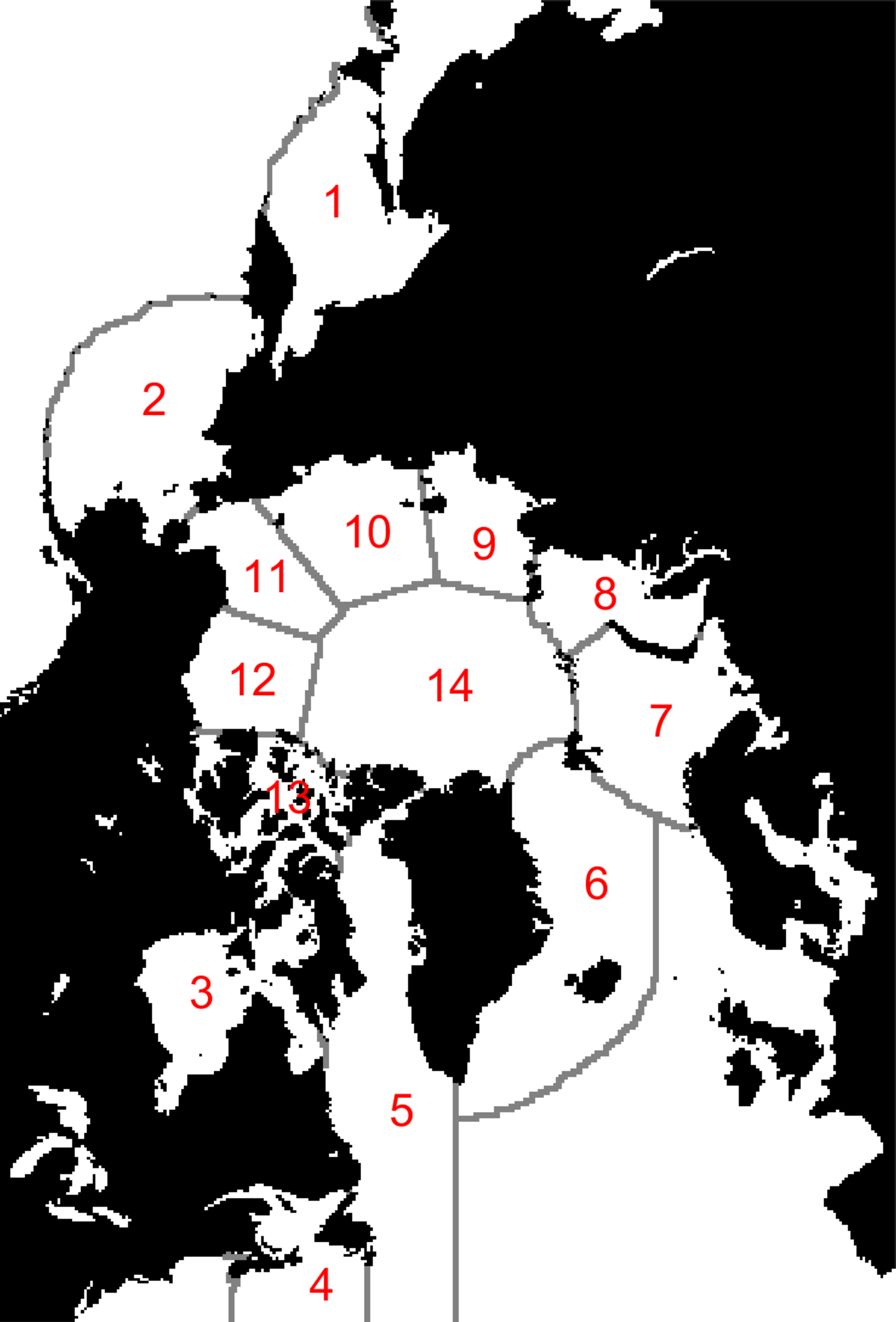}}
\caption{}
\end{subfigure}
\begin{subfigure}[b]{0.33\textwidth}
\centering
\setlength{\fboxsep}{0.9pt}
\fbox{\includegraphics[width=0.95\textwidth]{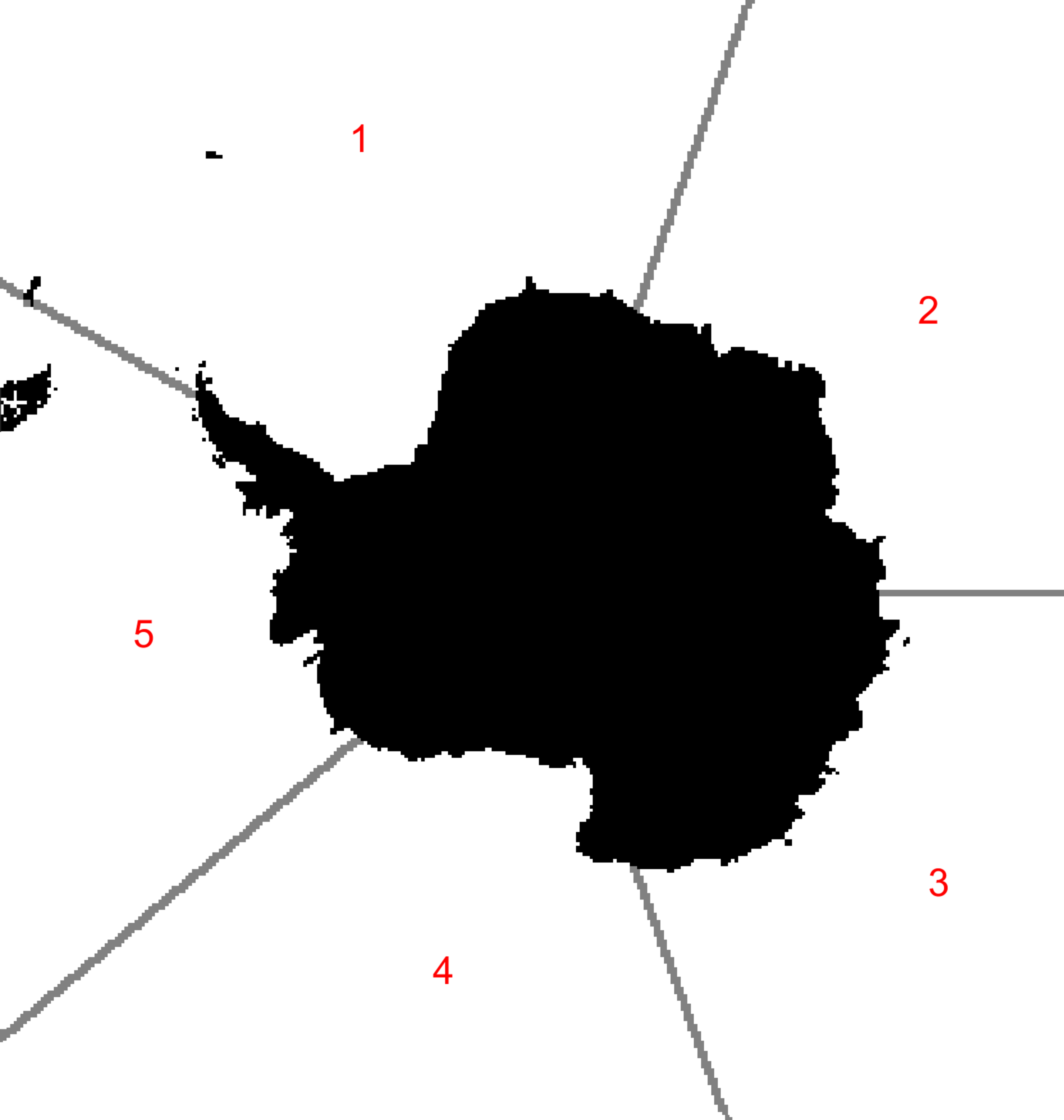}}
\caption{}
\end{subfigure}
\caption{Northern and southern hemisphere geographic regions considered. (a) Northern hemisphere geographic regions. 1: Sea of Okhotsk, 2: Bering Sea, 3: Hudson Bay, 4: Gulf of St Lawrence, 5: Baffin Bay, 6: East Greenland Sea, 7: Barents Sea, 8: Kara Sea, 9: Laptev Sea, 10: East Siberian Sea, 11: Chukchi Sea, 12: Beaufort Sea, 13: Canadian Arctic Archipelago, 14: Central Arctic Ocean. (b) Southern hemisphere geographic regions. 1: Weddell Sea, 2: Indian Ocean, 3: Pacific Ocean, 4: Ross Sea, 5: Bellingshausen and Amundsen Seas.}
\label{fig1_regions}
\end{figure}

A data set well suited for study and prediction with KMD is the satellite-based sea ice concentration measurements from the NSIDC Sea Ice Index \cite{nsidc_seaIceIndex}, due primarily to the long and continuous time period (from November 1978 to the present) and the large geographic regions over which this data is available. KMD analysis was applied to both entire Arctic and Antarctic sea region as a whole and to specific sea regions in each hemisphere.  The geographic regions used for the Arctic were those given by Boisvert and Stroeve \cite{boisvert2015arctic} and the Antarctic regions were those given by the NSIDC \cite{nsidc_seaIceIndex} (see Fig. \ref{fig1_regions} for the definitions of each region). Note that using these geographic definitions, not all of the ocean region in the northern hemisphere data is considered in the Arctic, therefore the sea ice concentration data from these non-Arctic regions were excluded from the KMD analysis.

Examination of a mode shows the geographic locations where the sea ice concentration has oscillatory, growth or decay behavior as determined by the associated eigenvalue, and thus allows one to associate particular Koopman modes with aspects of sea ice dynamics of interest. For example, modes with oscillatory frequencies with periods of one year correspond to the annual variation between the sea ice minimum extent in the late summer and the sea ice maximum extent in the late winter, and modes with eigenvalues near zero correspond to the mean sea ice concentration over the time period spanned by the data. Furthermore, the presence of eigenvalues with multi-year oscillatory periods suggests the presence of long-term periodic variations in sea ice behavior, and eigenvalues with real components leading to relatively slow growth or decay time constants suggests the existence of long-term increases or decreases in sea ice concentration, where again the associated modes show the geographic regions where such behaviors occur.  

Note that while, by definition, the modes resulting from KMD have complex-exponential time dependence, the identification modes possessing long term exponential decay is non-trivial because the identified modes have corresponding eigenvalues that are isolated in the eigenvalue space.  Linear growth or decay can be reproduced with a combination of exponential modes with eigenvalues close together in the eigenvalue space.  The long-term exponential decay modes identified here have no such clustering by the corresponding eigenvalues.

The Koopman modes and eigenvalues also permit prediction of the future sea ice concentration behavior. Data decomposed using KMD can be reconstructed over its original time period, and these same reconstruction equations allow for prediction of future behavior simply by increasing the time variable to values beyond the time period of the original input data. In this work we apply KMD to data over various multi-year time windows before a given year and produce prediction results for after the time window for which data was available, thus enabling a judgment of the goodness of the KMD predictions by comparison with the true “future” values. 

We applied KMD reconstruction techniques to the prediction of future sea ice concentrations in both the entire northern and southern hemispheres where data was available and also in particular geographic regions within each hemisphere. Because the dynamics of the seasonal variation in sea ice concentration differ greatly between high latitude regions, which have significant or complete sea ice coverage in the winter and can retain some sea ice through the summer, and lower latitude regions, which do not necessarily reach complete sea ice coverage in the winter nor retain any sea ice through the summer, it is of interest to examine the KMD prediction results in various seas, bays, and other specific regions. To that end, KMD prediction results for each of the sea regions described previously were examined separately.

The primary results of the work are:
\begin{itemize}
\item The presence of Koopman modes showing the change in geographic distribution of sea ice over the time period covered by satellite data, specifically the reduction in the mean Arctic and Antarctic sea ice concentration and the increased annual variation near West Antarctica and in the Arctic marginal seas.
\item Long-term exponential decay behavior in sea ice concentration in both the Arctic and Antarctic that is indicative of feedback mechanisms that accelerate decline in the extent of the sea ice cover \cite{comisoetal:2008,cohenetal:2014,goosseetal:2018}.
\item The ability of Koopman-based reconstruction techniques to predict future sea ice behavior over multi-year timescales.
\end{itemize}

\section*{Materials and Methods}

The objective of this study was to apply Koopman Mode Decomposition (KMD) and analysis techniques to existing satellite image data of sea ice concentration. Existing KMD algorithms were applied to the satellite images and the resulting Koopman modes and eigenvalues, defined in the following paragraph, showed the temporal and spatial details of the sea ice concentration dynamics.

The application of Koopman Mode Decomposition to time series data $g(t,\mathbf{z}_{0})\in \mathbb{R}^n\times I$, where $n$ is the number of field observations, $I$ the set of time snapshots and ${\bf z}_0$ an initial condition consists of expanding the time series observables  onto the Koopman eigenfunctions to produce a set of Koopman modes and Koopman eigenvalues \cite{Mezic:2005, Mezic:2013}:
\begin{equation*}
g(t,\mathbf{z}_{0})  = \sum_{j=1}^{\infty} e^{\lambda_{j} t} \mathbf{v}_{j}(\mathbf{z}_{0}) + e^{\bar{\lambda}_{j} t} \bar{\mathbf{v}}_{j}(\mathbf{z}_{0}) +{\bf n}(t)
\end{equation*}
where $\lambda_{j}$ are the Koopman eigenvalues, $\mathbf{v}_{j}$ are the Koopman modes, the overline indicates complex conjugation, and ${\bf n}(t)$ is part of the time evolution with continuous spectrum \cite{Mezic:2005,Mezic:2013}.  Note that the dependence on the initial state ${\bf z}_0$ is sometimes taken out of the mode itself, when eigenfunctions are used in the expansion. For each eigenvalue and its corresponding mode, the imaginary component of the eigenvalue $\Im(\lambda_{j}) = \omega_{j}$ determines the oscillation frequency $\omega_{j}$ of the mode $j$ (and the oscillatory period $\tau^{osc}_{j} = 1/\omega_{j}$) and the real component $\Re(\lambda_{j}) = 1/\tau^{decay}_{j}$ determines the growth or decay time constant $\tau^{decay}_{j}$ of the mode. It has been shown that short-timescale fluctuations can have a significant effect on large-scale features of climate systems, but the standard modeling practice of averaging the outputs of multiple simulation leads to the loss of this effect and thus incorrect predictions \cite{DingEtAl:2018}.  An advantage of KMD as a data-driven analytical technique is that in the eigenvalues it retains and uses the entire range of timescales present in the input data. The mode itself determines the spatial structure of the specific dynamical behavior given by the eigenvalue. 

The data pre-processing method was to convert the image data files from the NSIDC Sea Ice Index of sea ice concentration showing average monthly concentration to numerical arrays, remove the pixels corresponding to land areas, and reshape the remaining sea pixels into a 1-D array for each month. Note that there is a “polar data gap” in a circular region around the North Pole where concentration measurements are not available due to the coverage of the satellite based remote sensing instruments used to collect the sea ice concentration data. This region is traditionally either treated as completely ice covered or filled in based on the observed boundary conditions of the region \cite{StrongAndGolden:2016}. For KMD, it is not necessary to fill in this region and so the points in the polar data gap were excluded from our analysis.

Data files were missing for a small number (three in the northern hemisphere and two in the southern hemisphere) of months in the 1980's, so it was necessary to interpolate over the missing months to allow use of all data back to 1979, giving 40 full years of data (1979 to 2018).

The arrays for each month were then combined into a 2-D data matrix for performing KMD analysis run using KMD algorithms based on both Arnoldi \cite{Rowleyetal:2009} and DMD type methods \cite{DMDbook,Drmacetal:2018}. The results from the two categories of algorithms were found to be identical, which was taken to be a strong indication that the results are a good representation of the true Koopman eigenvalues and modes.

As described in the text, the calculated Koopman eigenvalues showed the time dependence (oscillatory and growth/decay) of the Koopman modes, which themselves showed the spatial structure of the time dependence of the input data.

To capture relatively short time scale dynamics, the analysis was performed on windowed data sets. The windowing consisted of performing KMD on subsets of the sea ice concentration data covering time periods of 5 to 40 years (e.g. five-year windows consisted of 1979-1983, 1980-1984, . . . , 2014-2018).

Reconstruction of the $N_{p}$ sea ice concentration pixel values $\mathbf{C}$ at discrete time step $k$ is performed using the Koopman eigenvalues $\lambda_{j}$ and the Koopman modes $\mathbf{v}_{j}$ obtained from applying KMD to the concentration values over $N$ time steps (months, in this case):
\begin{equation*}
\mathbf{C}_k = \sum_{j=1}^{N} \lambda_{j}^{k-1} \mathbf{v}_{j}
\end{equation*}
Here, there are $N$ Koopman eigenvalues and Koopman modes, where each Koopman eigenvalue is a single complex number and each Koopman mode has dimensions 1 by $N_{p}$.

For $1 \le k \le N$, $\mathbf{C_{k}}$ is termed a reconstruction of the $k$th time step in the original data $\mathbf{C}$, as the Koopman eigenvalues and modes came from a decomposition of the observations over this time range and should simply reproduce the data used as input to the KMD.  For $k>N$, $\mathbf{C_{k}}$ is a prediction of the future behavior of the sea ice concentration for the (future) $k$th time step, based on the system dynamics deduced from decomposition of earlier observations.

No probability distribution is assumed in the KMD process so no statistical methods were applied. The deviation between the KMD reconstruction-based predictions of future sea ice concentrations and the actual values are due to two factors: the finite dimensionality of numerical realizations of KMD algorithms, which for relatively high-dimensional data as used in this study is not expected to be a major source of error, and the stochastic nature of the underlying climatological processes driving sea ice concentration dynamics, which will produce behavior not predictable in a purely dynamical model such as that produced by KMD reconstruction-based predictions.
 
\section*{Results}

We have performed KMD processing on the sea ice concentration data using both Arnoldi-type and DMD-based KMD algorithms. All of the algorithms used gave very similar results. This suggests that the sea ice concentration data dynamical behavior is “well behaved” in the sense that the resulting condition number is sufficiently small that any of the various approximations of the Koopman decomposition are valid here \cite{Drmacetal:2018} and thus supports the conclusion that the KMD results obtained here are physically meaningful and not numerical artifacts.

\begin{figure}[h!]
\centering
\begin{subfigure}[b]{0.33\textwidth}
\includegraphics[width=0.96\textwidth]{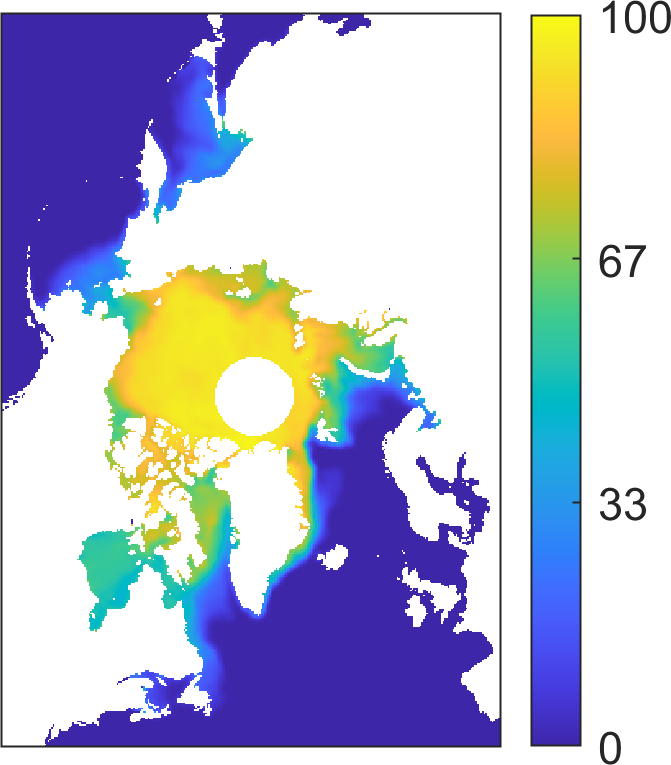}
\caption{}
\end{subfigure}
\begin{subfigure}[b]{0.33\textwidth}
\includegraphics[width=0.96\textwidth]{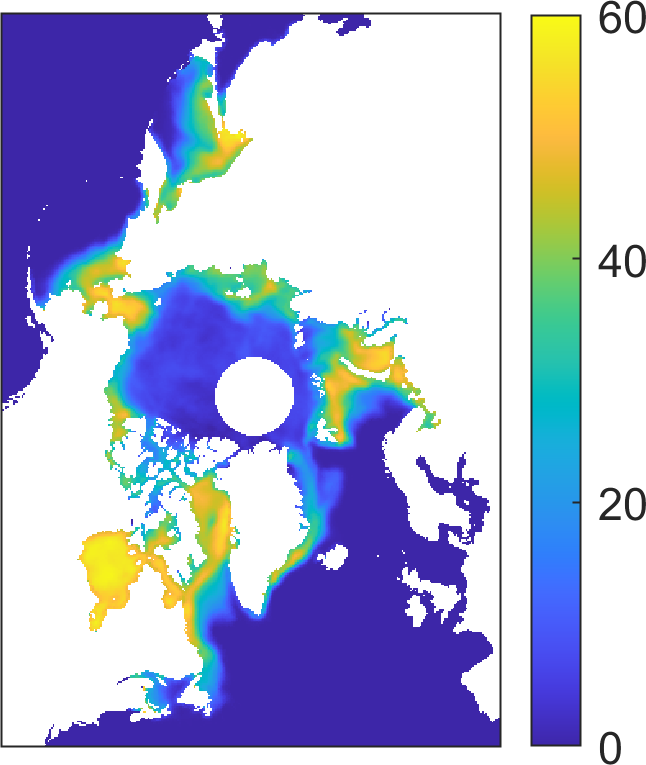}
\caption{}
\end{subfigure} \\
\centering
\begin{subfigure}[b]{0.33\textwidth}
\includegraphics[width=0.96\textwidth]{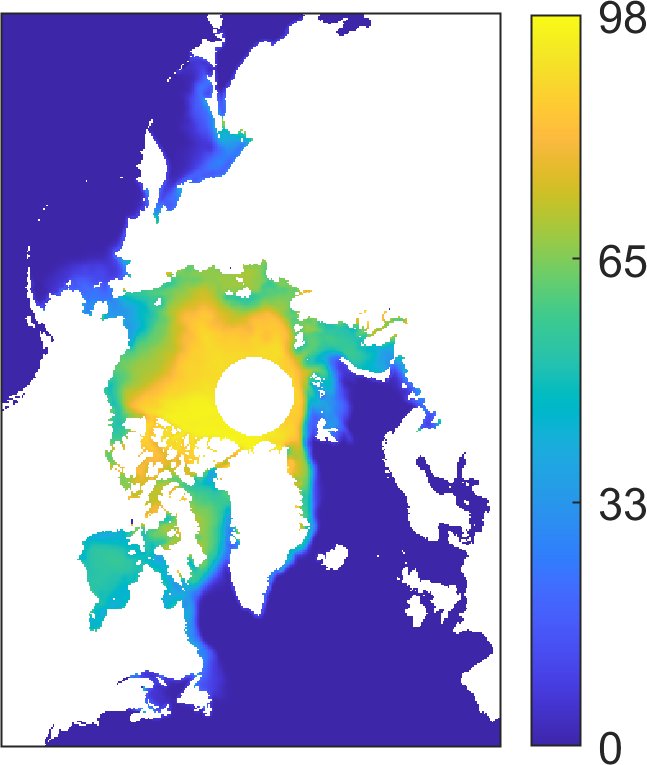}
\caption{}
\end{subfigure}
\begin{subfigure}[b]{0.33\textwidth}
\includegraphics[width=0.96\textwidth]{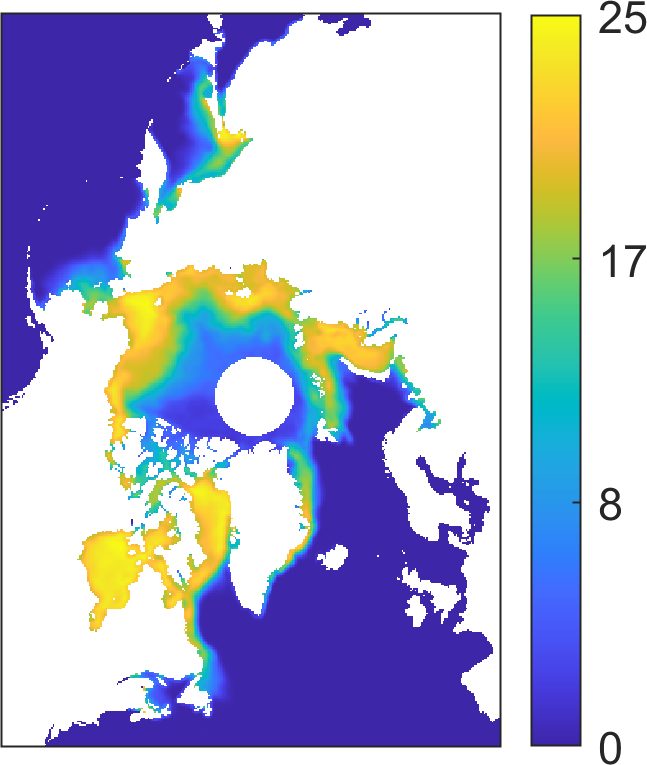}
\caption{}
\end{subfigure}
\caption{Koopman modes representing the mean and annual variation in sea ice
concentration over five year windows for the northern hemisphere. (a) Mean coverage,
1979-1983 period, (b) annual variation, 1979-1983 period, (c) mean coverage, 2014-2018 period, (d) annual variation, 2014-2018 period. The colorbar units are percent concentration.}
\label{modes_north}
\end{figure}

\begin{figure}[h!]
\centering
\begin{subfigure}[b]{0.33\textwidth}
\includegraphics[width=0.96\textwidth]{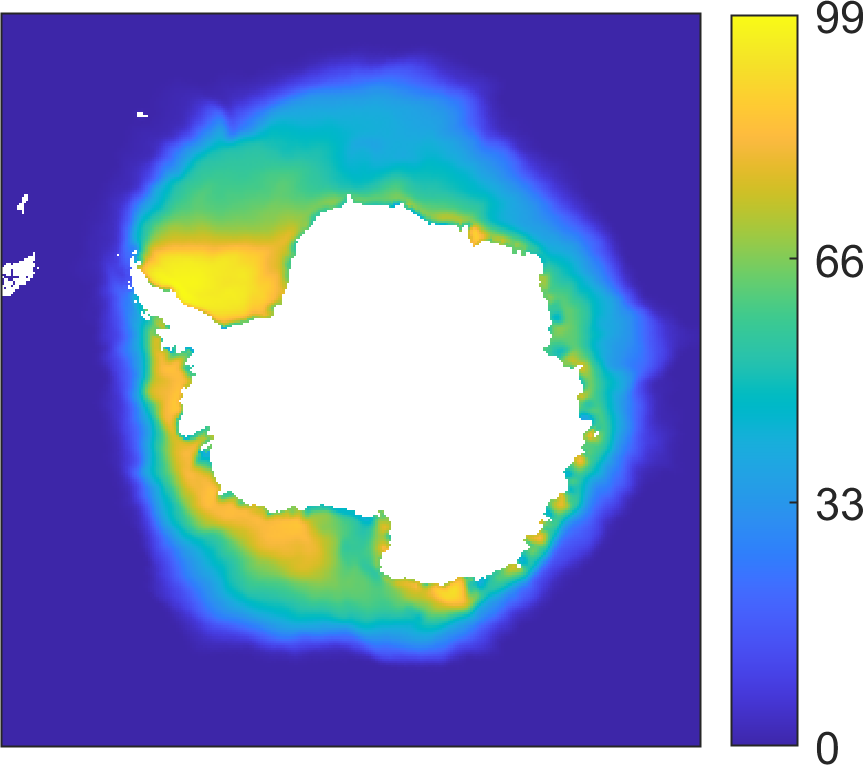}
\caption{}
\end{subfigure}
\begin{subfigure}[b]{0.33\textwidth}
\includegraphics[width=0.96\textwidth]{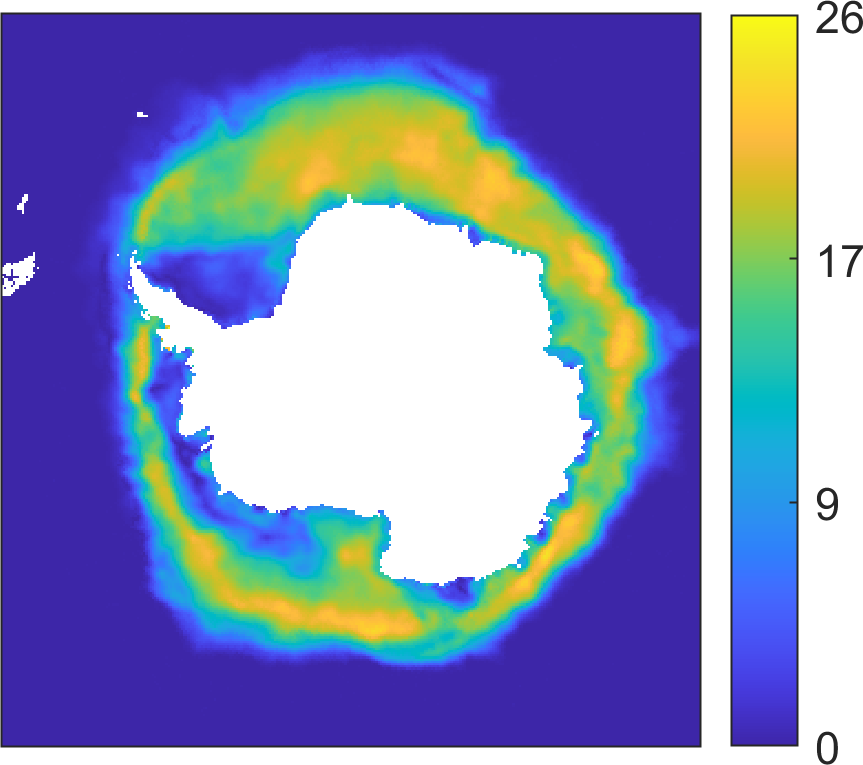}
\caption{}
\end{subfigure} \\
\centering
\begin{subfigure}[b]{0.33\textwidth}
\includegraphics[width=0.96\textwidth]{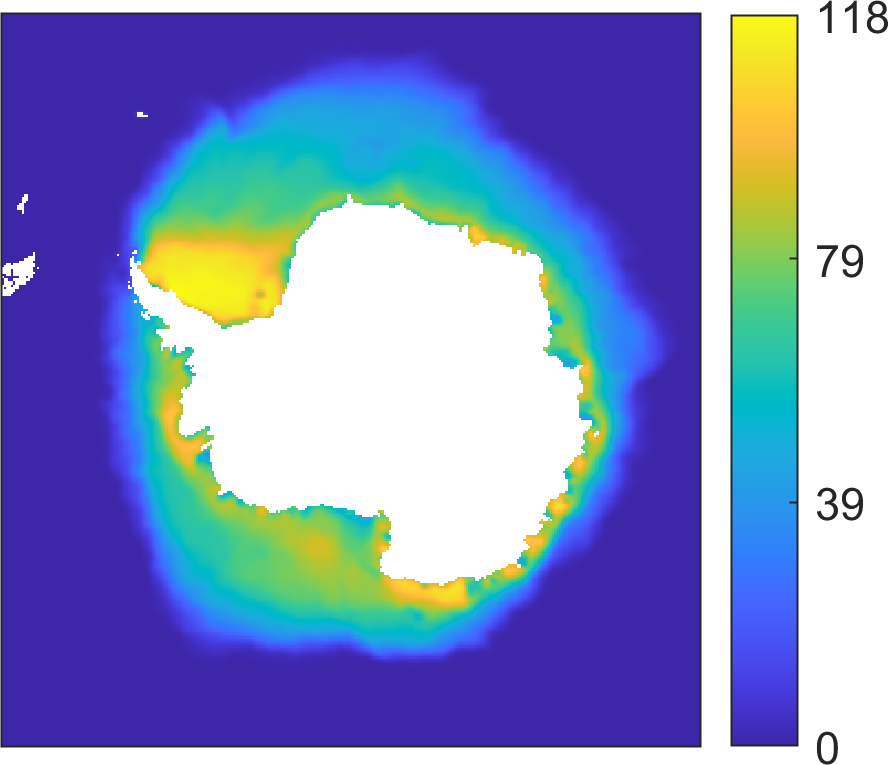}
\caption{}
\end{subfigure}
\begin{subfigure}[b]{0.33\textwidth}
\includegraphics[width=0.96\textwidth]{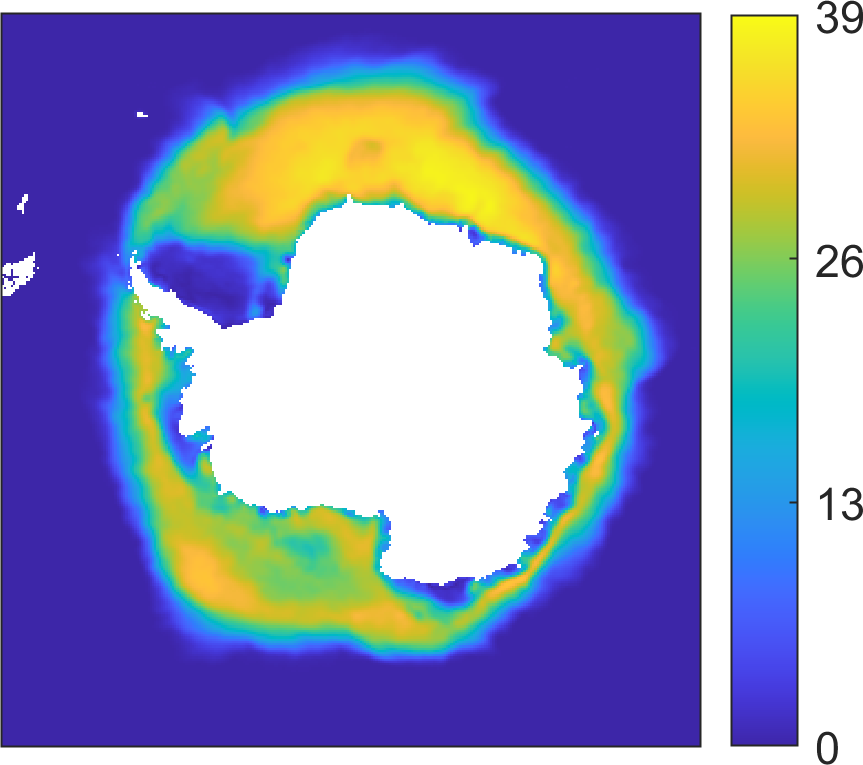}
\caption{}
\end{subfigure}
\caption{Koopman modes representing the mean and annual variation in sea ice
concentration over five year windows for the southern hemisphere. (a) Mean coverage,
1979-1983 period, (b) annual variation, 1979-1983 period, (c) mean coverage, 2014-2018 period, (d) annual variation, 2014-2018 period. The colorbar units are percent concentration, but as each mode is one component of a decomposition, concentration values in the modes can take on non-physical values.}
\label{modes_south}
\end{figure}

Figs. \ref{modes_north} and \ref{modes_south} show Koopman modes corresponding to the mean and annual variation in sea ice concentration for two 5-year time periods (1979-1983 and 2014-2018), as well as modes corresponding to long-term exponential decay over a 40-year time period. The mean mode in each case is defined as that with the largest L2-norm (taken over the components of the mode) and zero or negligible real and imaginary eigenvalue.  The annual mode in each case is the mode with a $\tau^{osc}$ value closest to 12 months. In all cases the annual mode was unambiguously identifiable as a large L2-norm mode.  The long-term exponential decay modes shown for the 40-year window are the two largest L2-norm modes with $\tau^{decay}$ periods greater than one year.  Note that although the modes have the same units as the input data, the modes can include non-physical values (i.e. concentration values less than 0\% or greater than 100\%) because the modes are mathematical structures resulting from a decomposition of the input observable data and not themselves direct representations of observable quantities.  A similar distinction can be made with, e.g., the Fourier coefficients in Fourier analysis, the values of which are not limited to values between the extrema of the analyzed function.

Comparison of the mean modes in Fig. \ref{modes_north} (a) with Fig. \ref{modes_north} (c), and Fig. \ref{modes_south} (a) with Fig. \ref{modes_south} (c), shows that the mean sea ice concentration is lower in the later time periods, suggesting that the sea ice does not reach as great an extent in the winter. Similarly, examination of the annual variation modes in the two time periods shows a greater annual variation in sea ice concentration in some regions.  This greater annual variation is particularly in the regions of the Beaufort Sea, Kara Sea, and the coastal corridor in between, and in the Bellingshausen and Amundsen Seas near West Antarctica and the Pacific Ocean in East Antarctica, with relatively little change in the Ross Sea area.  

The combination of the mean and annual variation modes shown in Fig. \ref{modes_north} and \ref{modes_south} can be viewed as first-order models of the sea ice concentration dynamics in each hemisphere over short annual timescales during the respective five-year windows. The observed decreases in sea ice concentration from the earlier to later periods suggests that a mode with primarily long-time decaying behavior is needed to reproduce the long-term loss of sea ice observed in the regions identified above. 

\begin{figure}[h!]
\centering
\begin{subfigure}[b]{0.33\textwidth}
\includegraphics[width=0.96\textwidth]{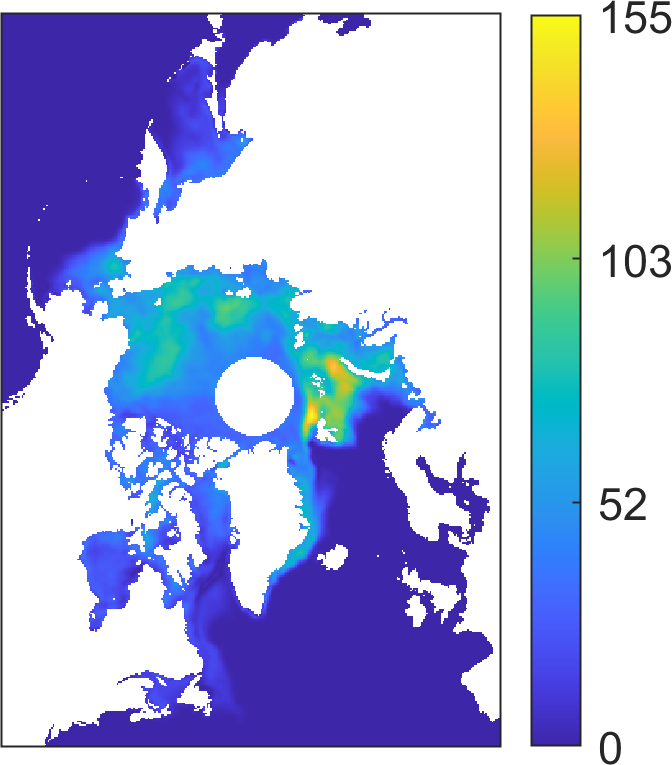}
\caption{}
\end{subfigure}
\begin{subfigure}[b]{0.33\textwidth}
\includegraphics[width=0.96\textwidth]{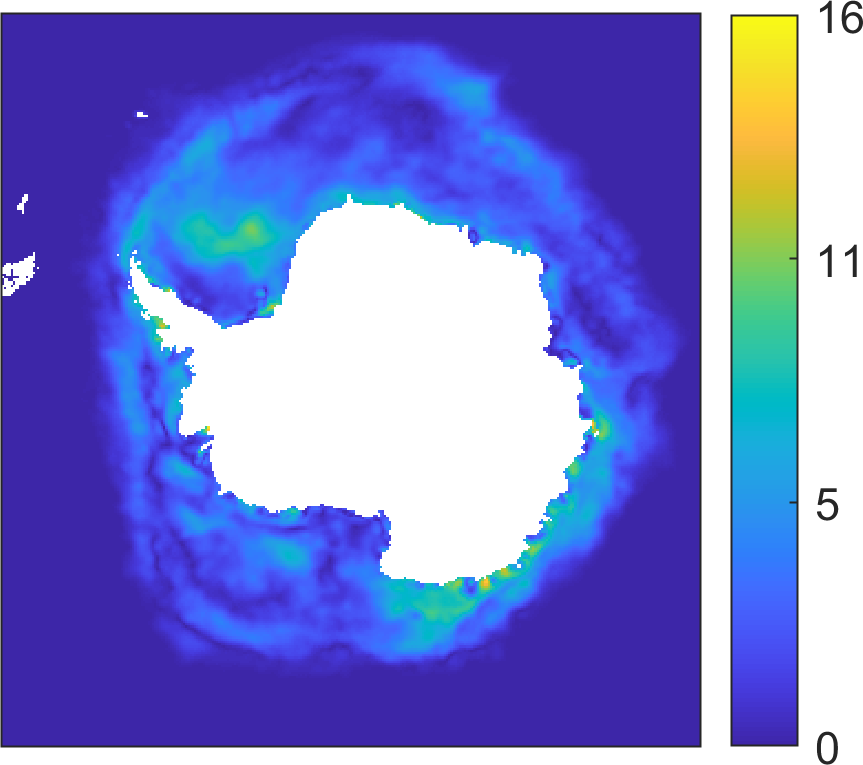}
\caption{}
\end{subfigure}
\caption{Koopman modes representing long term, exponential decay for the period 1979-2018. (a) Koopman mode representing long term, exponential decay in the northern hemisphere, corresponding to exponential decay with $\tau_{decay}$ =131 months. (b) Koopman mode representing long term, exponential decay in the southern hemisphere, corresponding to exponential decay in the Antarctic with $\tau_{decay}$ = 234 months. The colorbar units are percent concentration, but as each mode is one component of a decomposition, concentration values in the modes can take on non-physical values.}
\label{longterm}
\end{figure}

Such modes are apparent in analysis of the entire 40 year period of the data set. Fig. \ref{longterm} (a) shows such a mode from the northern hemisphere for the entire 40 year data set with $\tau_{decay}$ = 131 months and no oscillatory component. Consistent with other observations \cite{BeitschEtAl:2014, BoisvertEtAl:2016, RickerEtAl:2017}, these modes show that the decrease in sea ice coverage is most pronounced in the regions of the Beaufort Sea and the Arctic Ocean north of European Russia. This is also consistent with the changes in the mean and annual variation observed in the five-year window cases above. Fig. \ref{longterm} (b) shows an equivalent mode for the southern hemisphere, with slow decay ($\tau_{decay}$ = 234 months) and long oscillation period ($\tau^{osc}$ = 238 months) representing a decrease in sea ice concentration, primarily consisting of a decrease in sea ice concentration in West Antarctica. This region is known to be warming more rapidly than the region as a whole \cite{Rignot:2008},\cite{GardnerEtAl:2018}, so this KMD mode is consistent with that observation and the result from the five-year windows above showing decreased mean ice coverage and increased annual variation near West Antarctica. Note that the exponential decay described by the identified modes occurs in the geographic regions indicated by the spatial content of the mode, so the decrease in sea ice concentration on those time scales occurs locally in those regions.

The identification of an oscillatory period comparable to or greater than half of the time period of the input observations (that is, 238 months as described above, compared to the 480 months of available input data) does not violate the Nyquist criterion, as the large number of spatial dimensions effectively permits sampling of the underlying system dynamics over wider range of oscillation phase values than would observation of the time variation of a single spatial point.

\begin{figure}[h!]
\centering
\begin{subfigure}[b]{0.8\textwidth}
\includegraphics[width=\textwidth]{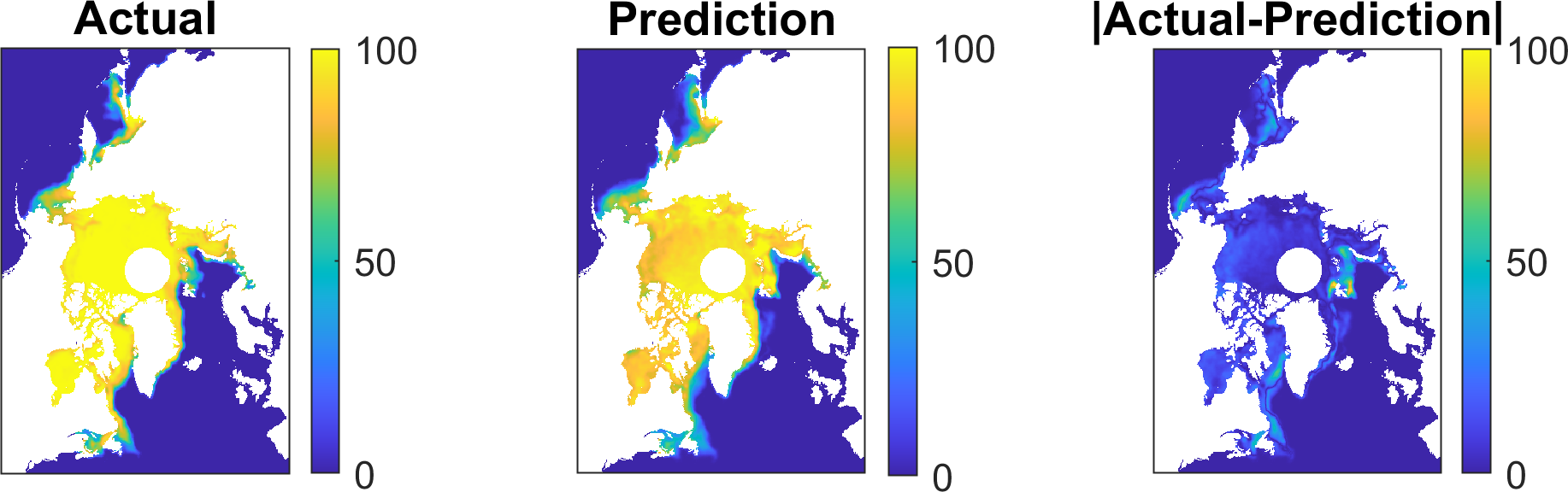}
\caption{}
\end{subfigure}
\begin{subfigure}[b]{0.8\textwidth}
\includegraphics[width=\textwidth]{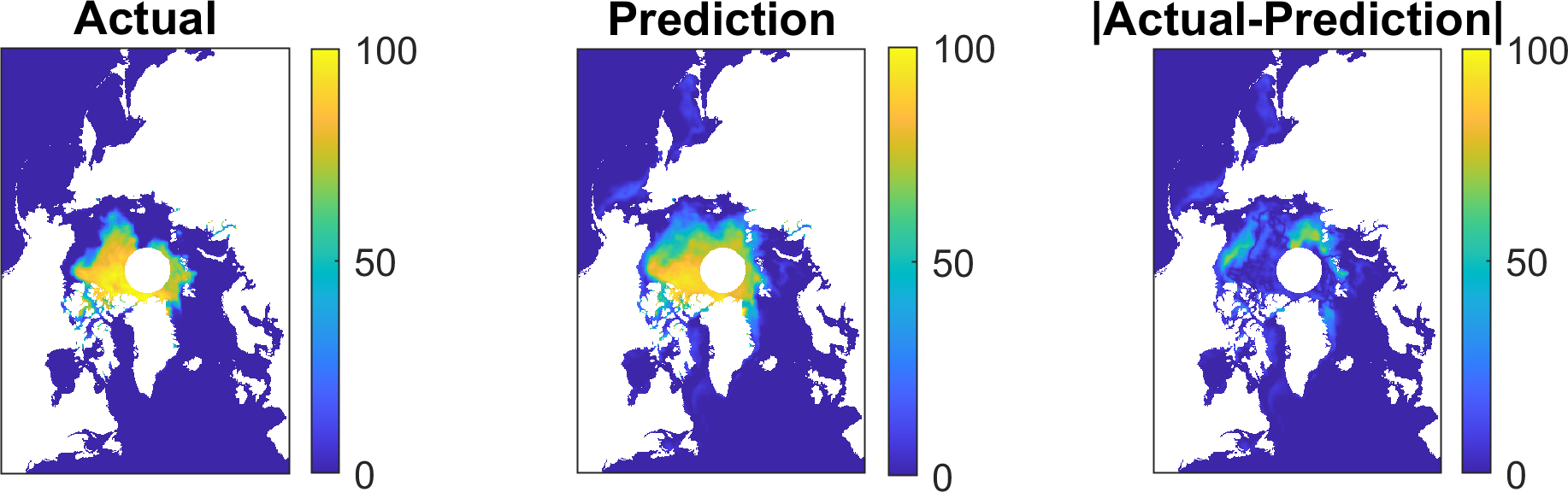}
\caption{}
\end{subfigure}
\begin{subfigure}[b]{0.8\textwidth}
\includegraphics[width=\textwidth]{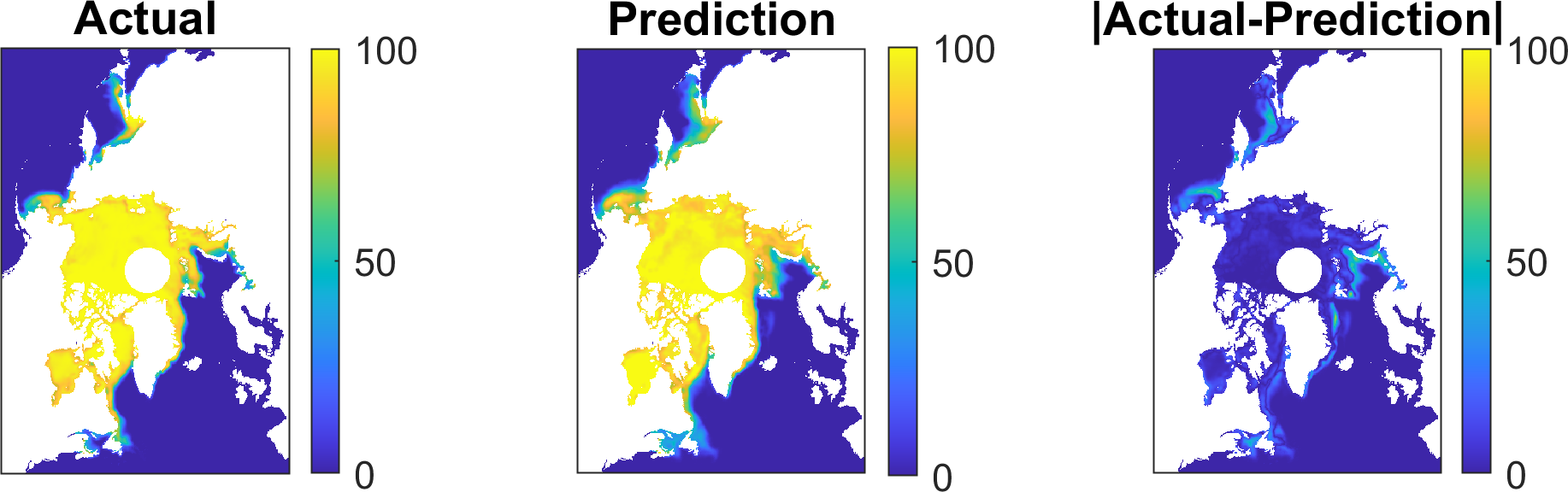}
\caption{}
\end{subfigure}
\begin{subfigure}[b]{0.8\textwidth}
\includegraphics[width=\textwidth]{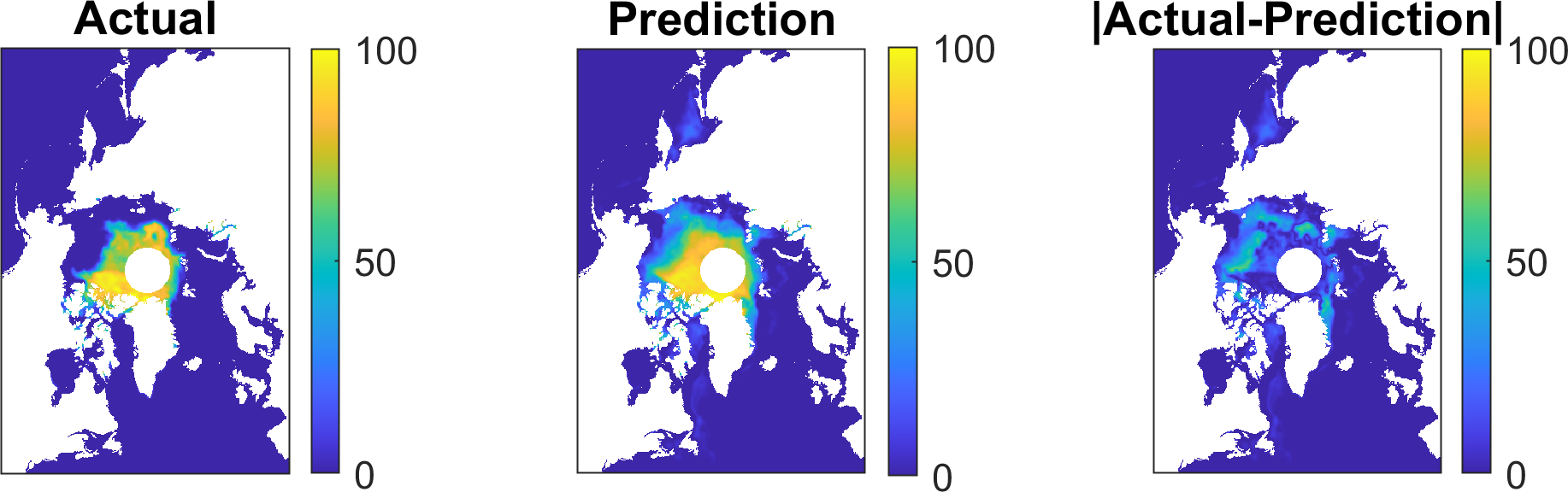}
\caption{}
\end{subfigure}
\caption{Comparison of actual data and prediction results for the winter sea ice maxima
and summer minima in the northern hemisphere for KMA performed on the input data period January 1984 to December 2013 (Left: actual concentration. Middle: predicted concentration. Right: Absolute difference between actual and predicted concentration). (a) March 2014, (b) September 2014, (c) March 2017, (d) September 2017.}
\label{prediction_north}
\end{figure}

Fig. \ref{prediction_north} shows example northern hemisphere sea ice concentration prediction results compared with the actual sea ice concentration for the same month. The data shown are the predictions for March, when the sea ice concentration is at its annual maximum, and September, when the sea ice concentration is at its annual minimum, for the four years following the 30-year window 1984-2013 that was used as the input data for KMD, where prediction calculated included all of the 360 KMD modes from the decomposition of the 360 months in the 30-year time window. It is seen that the prediction results match the general extent and magnitude of the actual data as well as capturing many small-scale features such as the shape of the concentration near the east coast of Greenland.

\begin{figure}[h!]
\centering
\begin{subfigure}[b]{0.85\textwidth}
\includegraphics[width=\textwidth]{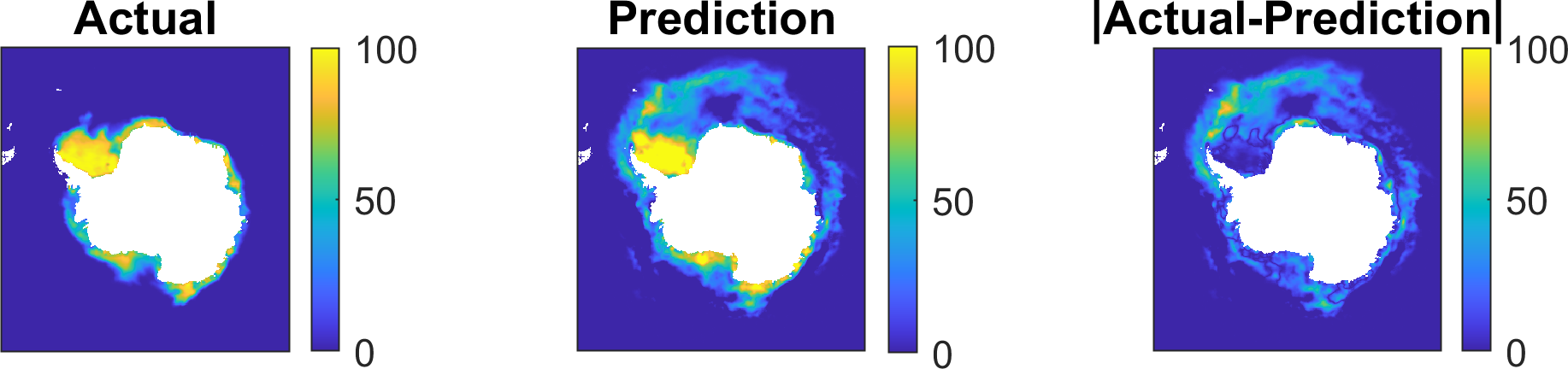}
\caption{}
\end{subfigure}
\begin{subfigure}[b]{0.85\textwidth}
\includegraphics[width=\textwidth]{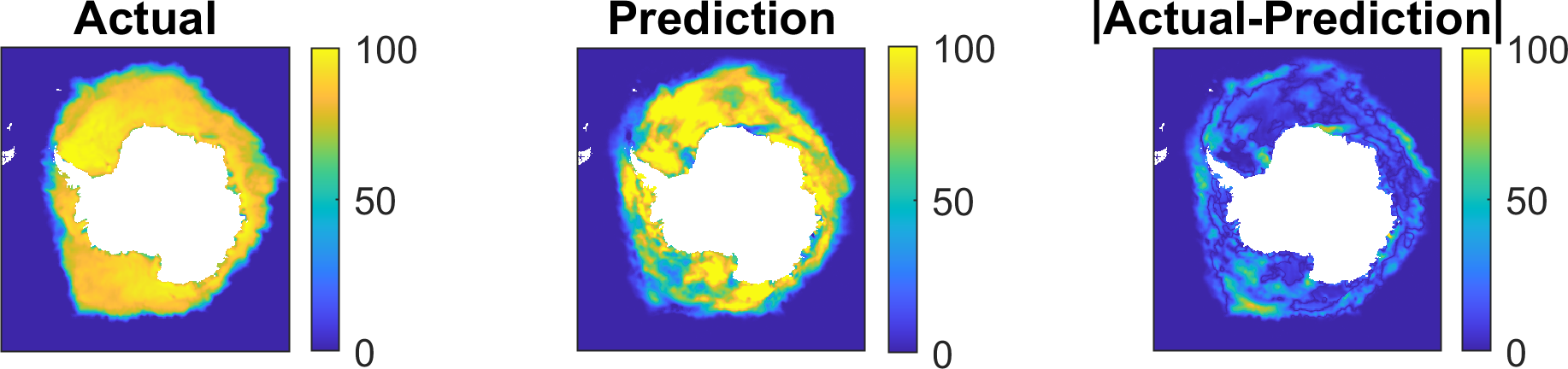}
\caption{}
\end{subfigure}
\begin{subfigure}[b]{0.85\textwidth}
\includegraphics[width=\textwidth]{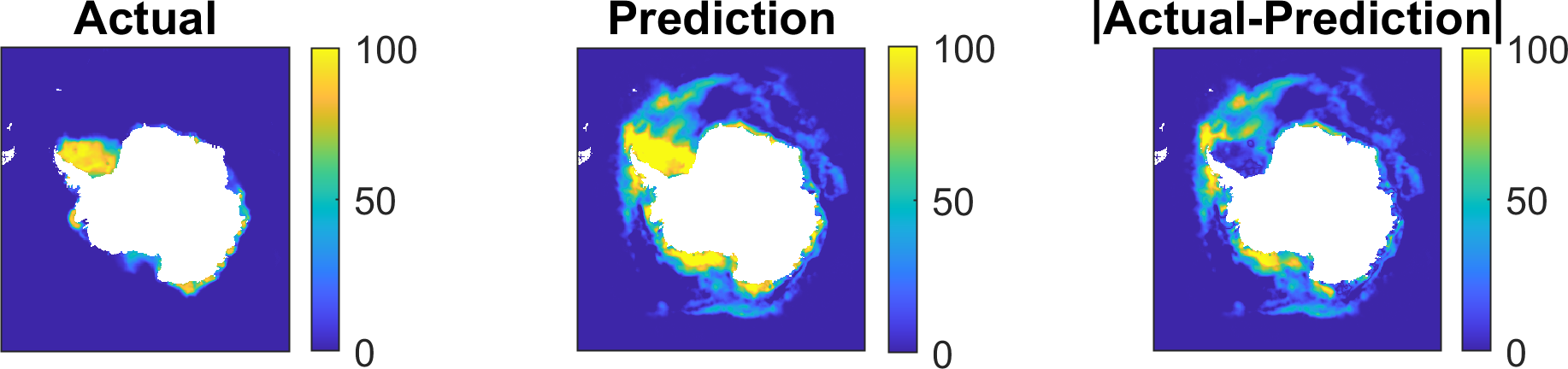}
\caption{}
\end{subfigure}
\begin{subfigure}[b]{0.85\textwidth}
\includegraphics[width=\textwidth]{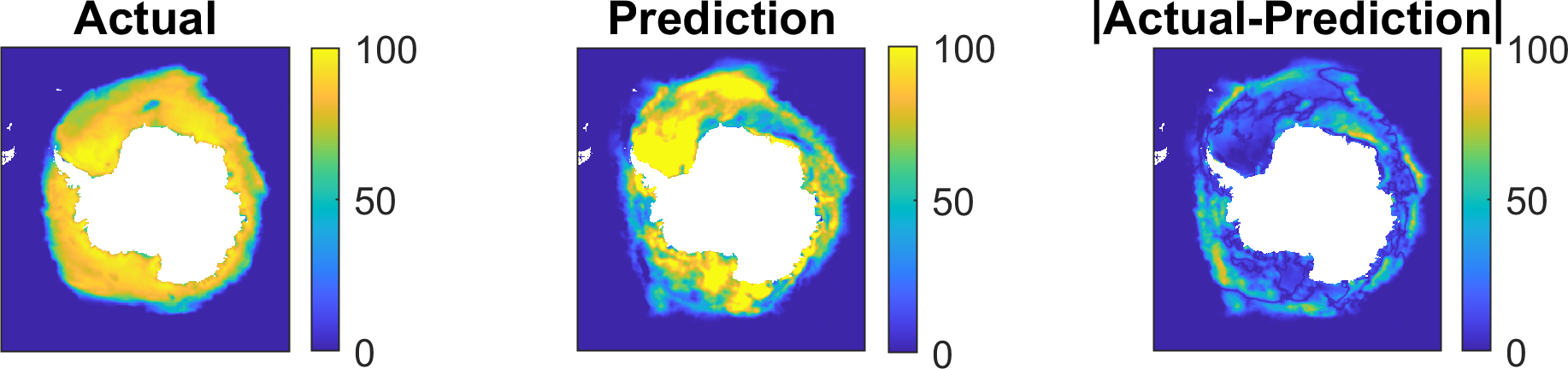}
\caption{}
\end{subfigure}
\caption{Comparison of actual data and prediction results for the summer sea ice minima and winter maxima in the southern hemisphere for 1984-2013 input data (Left: actual concentration. Middle: predicted concentration. Right: Absolute difference between actual and predicted concentration). (a) March 2014, (b) September 2014, (c) March 2017, (d) September 2017.}
\label{prediction_south}
\end{figure}

Fig. \ref{prediction_south} shows example southern hemisphere sea ice concentration prediction results compared with the actual sea ice concentration for the same month. The data shown are the predictions for March, when the sea ice concentration is at its annual minimum, and September, when the sea ice concentration is at its annual maximum, for the four years following the 30-year window 1984-2013 that was used as the input data for KMD. Again, it is seen that the prediction results match the general extent and magnitude of the actual data as well as capturing many small scale features, such as the remaining summer sea ice in the smaller bays and sea around East Antarctica and off of Marie Byrd Land in West Antarctica. The predictions tend to overestimate the sea ice concentration in the oceans away from the coast during the summer sea ice minima and underestimate the magnitude of the concentration values during the winter sea ice maxima while still capturing the winter sea ice extent well.

\begin{figure}[h!]
\centering
\begin{subfigure}[b]{0.65\textwidth}
\centering
\includegraphics[width=\textwidth]{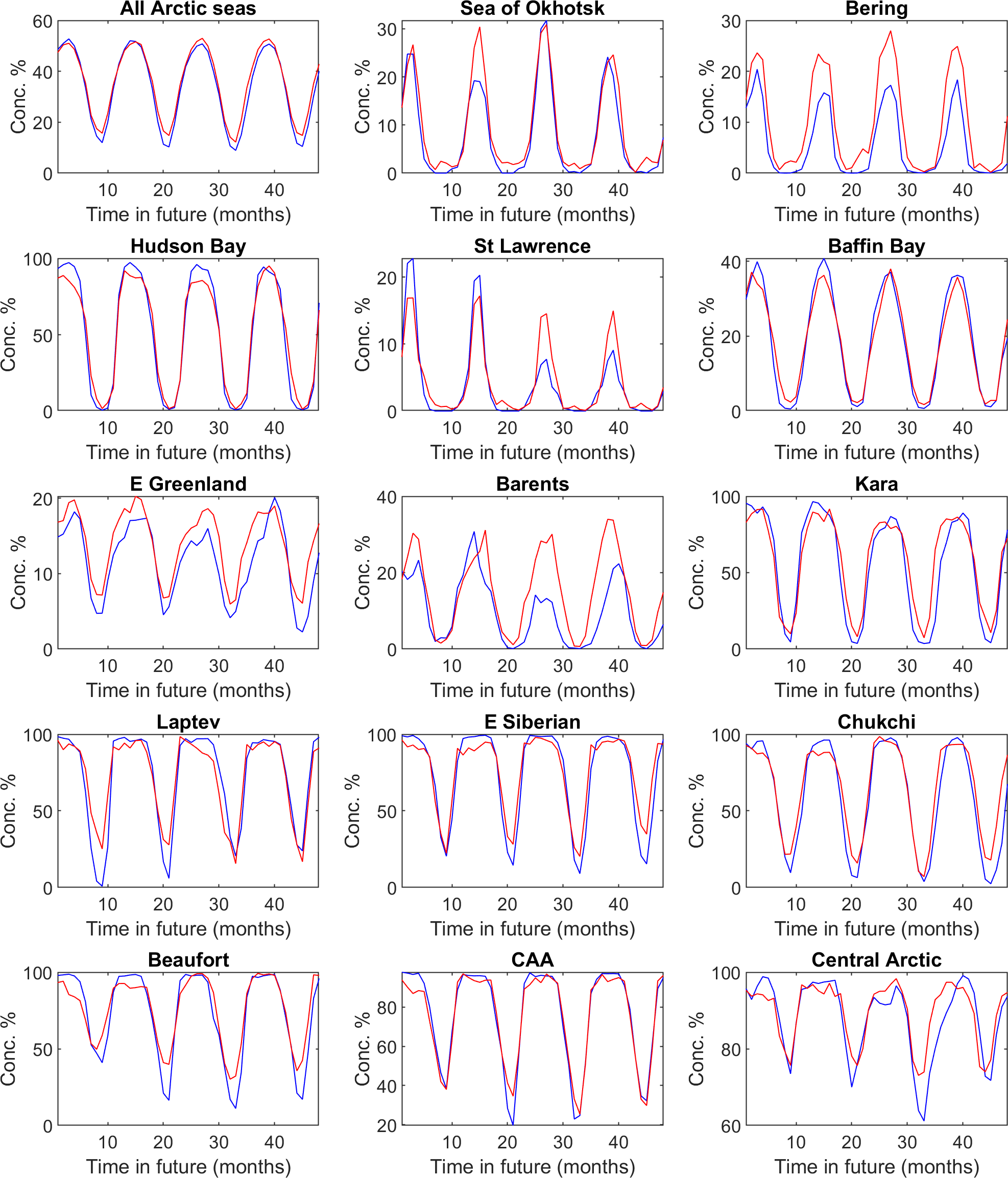}
\caption{}
\end{subfigure}\\
\begin{subfigure}[b]{0.65\textwidth}
\centering
\includegraphics[width=\textwidth]{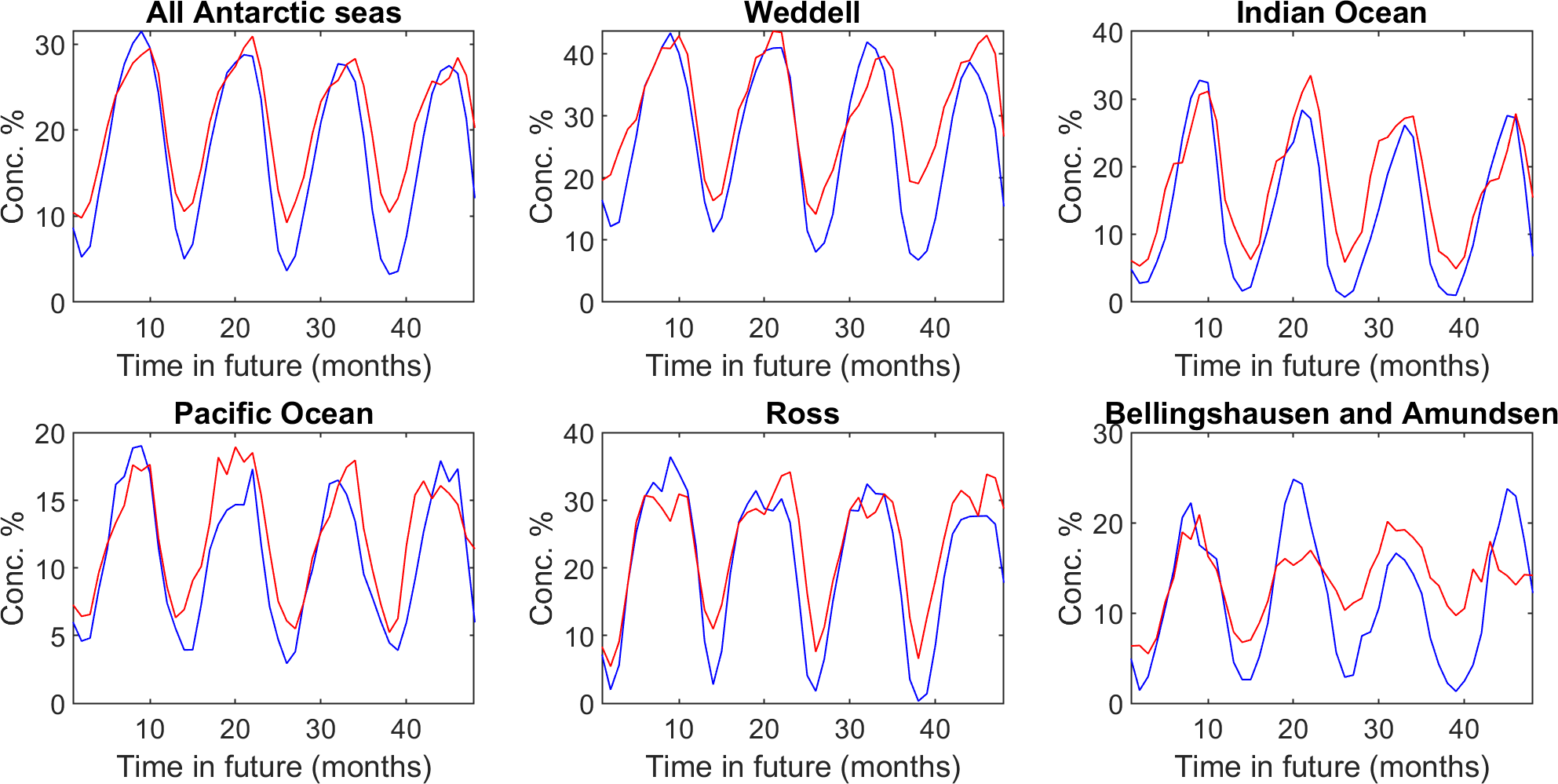}
\caption{}
\end{subfigure}
\caption{Prediction of mean sea ice concentration for all polar region seas. Comparison of actual (blue) and predicted (red) northern hemisphere mean sea ice concentration for the entire hemisphere and for each region.  Prediction is based on KMD of 30 year data range (1984 to 2013) and prediction of four future years.  Vertical axis units are sea ice concentration in percent. (a) Prediction for all Arctic polar region seas; (b) prediction for all Antarctic polar region seas.}
\label{prediction_mean}
\end{figure}

Fig. \ref{prediction_mean} shows a different view of the goodness of the northern hemisphere prediction, showing the mean of the actual data (blue lines) and predictions (red lines) for the entire northern hemisphere and for each region. The values were computed by averaging the actual or predicted sea ice concentration of each pixel within a given region. This shows the general trends of sea ice concentration in each region, rather than being a pixel-to-pixel comparison of the actual and predicted results. Here we see that the various Arctic seas mentioned above and other regions with large seasonal variations show good agreement between the actual and predicted results. This implies that while the prediction may not always be geographically precise in predicting the distribution of sea ice concentration within a particular region, it is generally successful at predicting the average sea ice concentration within the region. Close examination of regions of interest shows that even when discrepancies exist between the prediction and actual results of the summer sea ice concentration minimum, such as in the Central Arctic, the prediction does match the trend of the actual result (i.e., the summer minimum decreases year to year for the first three years, then increases in the fourth year). For the southern hemisphere it is seen that in the first year, the prediction of the maximum sea ice concentration is very good for each region, and even the minimum for the following year is reasonably well predicted in most of the regions.

\section*{Discussion and Conclusions}

These results show not only the previously-known existence of long term variation in sea ice concentration \cite{CavalieriAndParkinson:2012, StrongAndRigor:2013}, including long term decreases in sea ice coverage near West Antarctica and in the Arctic marginal seas, but also that a long-term exponential decrease in sea ice concentration exists and that Koopman Mode Decomposition allows a precise geographic view of where changes occur on annual and multi-year timescales. Such nonlinear trends in dynamics commonly result from positive feedback mechanisms such as those suggested in sea-ice dynamics studies \cite{comisoetal:2008,cohenetal:2014,goosseetal:2018}. The prediction results are possible over multi-year periods and also capture both large-scale features and trends in sea ice distribution and certain small-scale geographic details. The existence of long-term exponentially decaying modes seems to be of potentially substantial physical significance and warrants further measurement (including other physical fields) and investigation.

A limitation of the application of KMD is that, as a purely data-driven tool, it does not provide the physical insight into the underlying forcing or other drivers of a system’s dynamical behavior that numerical or theoretical models can provide. In this case, the geographic heterogeneity of the sea ice concentration behavior in Antarctica suggests a possible link with a proposed physical driver of the decrease in the Antarctic ice mass balance. Recent work \cite{rignot2019four} suggests the circumpolar deep water (CDW) flow as a physical mechanism for the decline of land-based ice due to increased glacier flow, where this increased decline is largest in the regions listed above where the decrease in the mean mode and increase in the annual variation mode is most apparent. The undersea topography of these regions is most consistent with the upwelling of relatively warm water by a strengthening CDW, leading to increase melting of ice shelves and, we suggest, reduced sea ice formation.

The accuracy of the predictions of future sea ice behavior by KMD reconstruction depends on the extent to which the sea ice concentration dynamics are governed by underlying nonlinear continuous processes, rather than stochastic or discontinuous drivers. That is, KMD-based prediction is expected to accurately predict variations in sea ice concentration due to the interactions of both long-term growth or decay and oscillatory behavior on fast or slow time-scales, however, changes due to “tipping points” such as the greater mixing between the Barents Sea and North Atlantic \cite{lind2018arctic} are not predictable from purely data-driven examination of sea ice concentration.

\section*{Data Management}

All data used in this work were obtained from the NSIDC Sea Ice Index \cite{nsidc_seaIceIndex}.

\section*{Acknowledgments}

This work was supported by ONR contracts N00014-18-P-2004 and N00014-19-C-1053 (Program Managers: Dr. Reza Malek-Madani and Dr. Behzad Kamgar-Parsi).

\nolinenumbers

\end{document}